\magnification=\magstep0
\documentstyle{amsppt}
\pagewidth{6.2in}
\input amstex
\topmatter
\title
Hypergeometric Series Acceleration Via the WZ method
\endtitle
\author
Tewodros Amdeberhan and Doron Zeilberger
\endauthor
\affil Department of Mathematics, Temple University,
Philadelphia PA 19122, USA \\
tewodros\@math.temple.edu, zeilberg\@math.temple.edu
\endaffil
\date Submitted: Sept 5, 1996. Accepted: Sept 12, 1996
\enddate
\dedicatory
Dedicated to Herb Wilf on his one million-first birthday
\enddedicatory
\abstract
Based on the WZ method, some series acceleration formulas are given. These formulas allow us to write down an infinite family of parametrized identities from any given identity of WZ type. Further, this family, in the case of the Zeta function, gives rise to many accelerated expressions for $\zeta(3)$.
\endabstract
\endtopmatter
\def\({\left(}
\def\){\right)}

\font\smcp=cmcsc8
\headline={\ifnum\pageno>1{\smcp the electronic journal of combinatorics 4 (2) (1997), \#R3\hfill\folio} \fi}
\document
{AMS Subject Classification:}
Primary 05A
\midspace{.05in}

We recall [Z] that a discrete function A(n,k) is called \it Hypergeometric \rm (or \it Closed Form \rm (CF)) \it in two variables \rm when the ratios $A(n+1,k)/A(n,k)$ and $A(n,k+1)/A(n,k)$ are both rational functions. A discrete 1-form $\omega=F(n,k)\delta k+G(n,k)\delta n$ is a \it WZ 1-form \rm if the pair (F,G) of CF functions satisfies $F(n+1,k) - F(n,k) = G(n,k+1) - G(n,k)$. \smallskip

We use: N and K for the \it forward shift operators \rm on n and k, respectively. $\Delta_n := N - 1$, $\Delta_k := K - 1$. \smallskip

Consider the WZ 1-form $\omega=F(n,k)\delta k+G(n,k)\delta n$. Then, we define the sequence $\omega_s, s=1,2,3,\dots$ of new WZ 1-forms: $\omega_s:=F_s\delta k+G_s\delta n$; where
$$F_s(n,k)=F(sn,k)\qquad\text{and}\qquad G_s(n,k)=\sum_{i=0}^{s-1}G(sn+i,k).$$

\bf Proposition: \rm  $\omega_s$ is WZ, for all s.\smallskip
\bf Proof: \rm  (a) $\omega_s$ is closed:
$$\align
\Delta_nF_s&=F(s(n+1),k)-F(sn,k)\\
&=\sum_{i=0}^{s-1}\biggl(F(sn+i+1,k)-F(sn+i,k)\biggr)\\
&=\sum_{i=0}^{s-1}\biggl(G(sn+i,k+1)-G(sn+i,k)\biggr)\\
&=\sum_{i=0}^{s-1}G(sn+i,k+1)-\sum_{i=0}^{s-1}G(sn+i,k)\\
&=\Delta_kG_s.
\endalign$$

\pagebreak
\midspace{.4in}

Note that since $\omega$ is a WZ, it has the form ([Z], p.590): 
$$\omega=f(n,k)\bigl(P(n,k)\delta k+Q(n,k)\delta n\bigr)\tag{$*$}$$
for some CF f and some polynomials P and Q. \smallskip

(b) $\omega_s$ has the form $(*)$:\par

Indeed, $\omega_s$ can be rewritten as: 
$$\align
\omega_s&=f(sn,k)\biggl(P(sn,k)\delta k+\sum_{i=0}^{s-1}\frac{f(sn+i,k)}{f(sn,k)}Q(sn+i,k)\delta n\biggr)\\
&=f(sn,k)\bigl(P(sn,k)\delta k+R(n,k)\delta n\bigr);
\endalign$$
where R(n,k) is a rational function and f(sn,k) is still CF. Hence after pulling out a common denominator, we see that $\omega_s$ too has the form $(*)$. This proves the Proposition. $\square$\bigskip

\bf Theorem 1: \rm ([Z], Theorem 7, p.596) For any WZ pair (F,G)
$$\sum_{n=0}^{\infty}G(n,0) = \sum_{n=1}^{\infty}\(F(n,n-1)+G(n-1,n-1)\)
-\lim_{n\to\infty}\sum_{k=0}^{n-1}F(n,k)
,$$
whenever both side converge.\smallskip

\bf Formula 1:\rm 
$$\sum_{n=0}^{\infty}G(n,0)=\sum_{n=0}^{\infty}\biggl(F(s(n+1),n)+\sum_{i=0}^{s-1}G(sn+i,n)\biggr)
-\lim_{n\to\infty}\sum_{k=0}^{n-1}F(sn,k)
.\tag 1$$

\bf Proof: \rm Apply Theorem 1 above on $\omega_s$. Alternatively, integrate $\omega$ along the boundary contour $\partial\Omega_s$ of the region $\Omega_s=\{(n,k): sn\geq k\}.$ $\square$\bigskip

\bf Formula 2: \rm We also have that
$$\sum_{k=0}^{\infty}F(0,k)-\lim_{n\to\infty}\sum_{k=0}^{n}F(n,k) = \sum_{n=0}^{\infty}G(n,0)-\lim_{k\to\infty}\sum_{n=0}^{k}G(n,k),\tag2$$
whenever both side converge.\smallskip

\bf Proof: \rm Integrate $\omega$ along the boundary contour $\partial\Omega_0$ of the region $\Omega_0=\{(n,k): n\geq 0, k\geq 0\}.$ $\square$\smallskip

\pagebreak
\midspace{.4in}

\bf Remark: \rm By shear symmetry, a formulation similar to (2) can be given in `k'. And a combination leads to: \smallskip

\bf Formula 3: \rm For $\omega_{s,t}=F_{s,t}\delta k+G_{s,t}\delta n$; where $$F_{s,t}(n,k)=\sum_{j=0}^{t-1}F(sn,tk+j)\qquad\text{and}\qquad G_{s,t}(n,k)=\sum_{i=0}^{s-1}G(sn+i,tk),\qquad\text{we have}$$
$$\sum_{n=0}^{\infty}G(n,0)=\sum_{n=0}^{\infty}\biggl(\sum_{j=0}^{t-1}F(s(n+1),tn+j)+\sum_{i=0}^{s-1}G(sn+i,tn)\biggr)
-\lim_{n\to\infty}\sum_{k=0}^{n-1}F_{s,t}(n,k).\tag3$$
\bigskip
 
Analogous statements hold in several variables. To wit: \smallskip

for the WZ 1-form in 3 variables, $\omega_{s,t,r}:=F_{s,t,r}\delta k+G_{s,t,r}\delta n+H_{s,t,r}\delta a$; where
$$\alignat3
F_{s,t,r}(n,k,a)&=\sum_{j=0}^{t-1}F(sn,tk+j,ra),&\qquad G_{s,t,r}(n,k,a)&=\sum_{i=0}^{s-1}G(sn+i,tk,ra)\qquad\text{and}\\
H_{s,t,r}(n,k,a)&=\sum_{u=0}^{r-1}H(sn,tk,ra+u),\endalignat$$

\bf Formula 4: \rm
$$
\sum_{n=0}^{\infty}H(0,0,n)=\sum_{n=0}^{\infty}\biggl(\sum_{u=0}^{r-1}H(s(n+1),t(n+1),rn+u)+\sum_{j=0}^{t-1}F(s(n+1),tn+j,rn)+\sum_{i=0}^{s-1}G(sn+i,tn,rn)\biggr)$$
$$-\lim_{a\to\infty}\sum_{k=0}^{a+1}F_{s,t,r}(a+1,k,a)
-\lim_{a\to\infty}\sum_{n=0}^{a+1}G_{s,t,r}(n,a+1,a).$$
\smallskip
In [A], formula (1) was used to give a list of series acceleration for $\zeta(3)$ (where $F(n,k)$ is given and its companion G(n,k) is produced by the amazing Maple Package {\tt EKHAD} accompanying [PWZ]). A small Maple Package {\tt accel} applying (3) is available at 
{\tt http://www.math.temple.edu/\~{}[tewodros, zeilberg]}.\smallskip

For example: with $F(n,k) = (-1)^k \frac{n!^6(2n-k-1)!k!^3}{2(n+k+1)!^2(2n)!^3}$, s=1 and t=1 {\tt accel} produces the following pretty formula:
$$\zeta(3)=\sum_{n=0}^{\infty}(-1)^n\frac{n!^{10}(205n^2+250n+77)}{64(2n+1)!^5}.\tag{$**$}$$
Greg Fee and Simon Plouffe used $(**)$ in their evaluation of $\zeta(3)$ to \it 520,000 digits \rm (available at {\tt http://www.cecm.sfu.ca/projects/ISC/records.html}).
\midspace{.05in}

\bf ACKNOWLEDGMENT: \rm We would like to express our gratitude to Professor Herbert Wilf for his valuable comments and suggestions.

\pagebreak
\midspace{.4in}

\enddocument